\documentstyle[12pt]{amsart}

\theoremstyle{plain}
  \newtheorem{thm}{Theorem}[section]
  \newtheorem{cor}[thm]{Corollary}
  \newtheorem{lem}[thm]{Lemma}
  \newtheorem{prop}[thm]{Proposition}
\theoremstyle{definition}
  \newtheorem{defn}{Definition}[section]
\theoremstyle{remark}
  \newtheorem{rem}{Remark}[section]
  \newtheorem{ack}{Acknowledgement}

\numberwithin{equation}{section}
\def\a{\alpha}
\def\b{\beta}
\def\d{\delta}
\def\f{\phi}
\def\g{\gamma}
\def\h{\eta}
\def\k{\kappa}
\def\l{\lambda}
\def\m{\mu}
\def\n{\nu}
\def\p{\pi}
\def\t{\tau}
\def\x{\xi}
\def\ve{\varepsilon}
\def\vt{\vartheta}
\def\vt1{\vartheta_1}
\def\L{\Lambda}
\def\otimeshat{\overset{\wedge}{\otimes}}
\def\Rcheck{\overset{\vee}{R}}
\def\Wcheck{\overset{\vee}{W}}
\def\pbar{\bar{\f}}
\def\rt-1{\sqrt{-1}}
\def\bk{\bold k}
\def\pk{\f_{\bk}}
\def\pkbar{\bar{\f}_{\bk}}
\def\Lcheckk{\overset{\vee}{L}_{\bk}}
\def\V{{\cal V}}
\def\Vm{{\cal V^-}}
\def\W{{\cal W}}
\def\tVk{{\tilde{V}_k}}
\def\vpk{{\varphi_{\bk}}}
\def\vpkbar{{\bar{\varphi}_{\bk}}}
\begin{document}
\title[Vertex--IRF correspondence for an elliptic $R$-operator]
{Vertex--IRF correspondence and factorized\\
L-operators for an elliptic ${\bold R}$-operator}
\author{Youichi SHIBUKAWA}
\address{Department of Mathematics,
Hokkaido University,
Sapporo 060 JAPAN}
\email{shibu@@math.hokudai.ac.jp}
\date{September 8, 1994\\
Revised December 26, 1994}
\maketitle
\begin{abstract}
  As for an elliptic $R$-operator which satisfies the
  Yang--Baxter equation, the incoming
  and outgoing intertwining vectors are constructed, and
the vertex--IRF correspondence for the elliptic $R$-operator
is obtained.
The vertex--IRF correspondence implies that the Boltzmann
weights of the IRF
model satisfy the star--triangle relation.
By means of these intertwining vectors, the
  factorized L-operators for the elliptic $R$-operator
are also constructed.
The vertex--IRF correspondence and the factorized L-operators
for Belavin's $R$-matrix are reproduced from those of the
elliptic $R$-operator.
\end{abstract}
\setcounter{section}{-1}
\section{Introduction}
\footnote[0]{1991 {\it Mathematics Subject
    Classification\/}. 17B37, 81R50, 82B23.
\endgraf
{\it Key words and phrases\/}. elliptic $R$-operator,
vertex--IRF correspondence,
  intertwining vector, factorized L-operator,
  Belavin's $R$-matrix.}
In~\cite{shibukawa1,shibukawa2,shibukawa3}
we have introduced an infinite-dimensional $R$-matrix.
It is a new solution of the
Yang--Baxter equation.
By means of the Fourier transformation of the $R$-matrix, we
defined an $R$-operator acting on some function space. This
$R$-operator also satisfies the Yang--Baxter equation.
Since this operator is deeply linked to analytic properties
of an elliptic theta
function, we call it the elliptic $R$-operator.
We have shown some properties satisfied by the elliptic
$R$-operator, for example, first inversion relation, fusion
procedure, etc.
For the trigonometric degenerate case of the elliptic
$R$-operator, we proved that the finite-dimensional,
trigonometric $R$-matrices
are constructed from the $R$-operator through restricting
the domain of the $R$-operator to some finite-dimensional subspaces.
Recently Felder and Pasquier~\cite{felder} showed that
Belavin's $R$-matrix \cite{belavin,richey} can be obtained through
restricting the domain of a
modified version of the
elliptic $R$-operator to a suitable finite-dimensional
subspace.

In \cite{baxter}, Baxter has introduced the intertwining
vectors for the eight--vertex model.
Jimbo, Miwa and Okado~\cite{jimbo}
constructed the outgoing intertwining vectors between
Belavin's vertex model and the $A_{n-1}^{(1)}$ face model.
We call this relation the vertex--IRF correspondence for
Belavin's $R$-matrix.
Hasegawa~\cite{hasegawa1,hasegawa2}, Quano and
Fujii~\cite{quano}
defined the incoming intertwining vectors which are the
dual vectors of the outgoing intertwining vectors.
Then they constructed the factorized L-operators for
Belavin's $R$-matrix.
The vertex--IRF correspondence plays a central role in their
methods.

The aim of this paper is to extend the result above to the
elliptic $R$-operator.

Our strategy to construct factorized L-operators for the
elliptic $R$-operator is as
follows.
At first we define incoming intertwining vectors
$\pbar_\l^\k$ of the elliptic $R$-operator $\Rcheck(\x)$ and
establish a vertex--IRF correspondence.
The vertex--IRF correspondence plays the most important role
in this paper.
Next we find finite-dimensional subspaces with the following
property (cf.~Theorem
\ref{thm:felder});
\[
\Rcheck(\x_{12})(V_k(\x_1)\otimes V_k(\x_2+\m))
\subset V_k(\x_2)\otimes V_k(\x_1+\m),
\]
where $\x_{12}:=\x_1-\x_2$.
Then we
define outgoing intertwining
vectors $\pk(\x)_\l^\k(z) \in V_k(\x+|\l|_{\bk})$, which are the duals
of $\pbar_\l^\k \Big\vert_{V_k(\x+|\l|_{\bk})}$.
Making use of the properties of the incoming and outgoing
intertwining vectors, we can easily construct factorized
L-operators.

This paper is organized as follows.
In Section 1, we review the properties of
the elliptic
$R$-operator $\Rcheck(\x)$ proved in
\cite{shibukawa1,shibukawa2,shibukawa3,felder}.
In Section 2, we shall
define incoming intertwining vectors $\pbar_\l^\k$
and Boltzmann
weights of an IRF model
\(\Wcheck
\begin{bmatrix}
  &\k'&\\
\l&\x&\n\\
&\k&
\end{bmatrix}
\).
Then we have the vertex--IRF correspondence for the elliptic
$R$-operator (Theorem~\ref{thm:vertex-irf}).
\begin{thm}[Vertex--IRF correspondence]
For $\l,\k,\n\in\L$
\[ \pbar_\l^\k \otimes \pbar_\k^\n \Rcheck(\x)
=\sum_{\k'\in\L}
\Wcheck
\begin{bmatrix}
  &\k'&\\
\l&\x&\n\\
&\k&
\end{bmatrix}
\pbar_\l^{\k'}\otimes \pbar_{\k'}^\n \ .
\]
\end{thm}
{}From this theorem, we can show that the Boltzmann weights
satisfy the star--triangle relation
(Corollary~\ref{cor:star-triangle}).
\begin{cor}
  The Boltzmann weights of the IRF model satisfy the
  star--triangle relation\rom;
  \begin{equation*}
    \begin{split}
      \sum_{\k'\in\L}
&
\Wcheck
\begin{bmatrix}
  &\k'&\\
\k&\x_{12}&\g\\
&\n&
\end{bmatrix}
\Wcheck
\begin{bmatrix}
  &\a&\\
\l&\x_{13}&\k'\\
&\k&
\end{bmatrix}
\Wcheck
\begin{bmatrix}
  &\b&\\
\a&\x_{23}&\g\\
&\k'&
\end{bmatrix}
\\
& \qquad =
\sum_{\k'\in\L}
\Wcheck
\begin{bmatrix}
  &\k'&\\
\l&\x_{23}&\n\\
&\k&
\end{bmatrix}
\Wcheck
\begin{bmatrix}
  &\b&\\
\k'&\x_{13}&\g\\
&\n&
\end{bmatrix}
\Wcheck
\begin{bmatrix}
  &\a&\\
\l&\x_{12}&\b\\
&\k'&
\end{bmatrix}
,
    \end{split}
  \end{equation*}
for $\l,\k,\n,\a,\b,\g\in\L$.
\end{cor}
This IRF model can be regarded as the limiting case $n\to
\infty$ of the $A_{n-1}^{(1)}$ face model.
In Section 3, making use of the results obtained by Felder
and Pasquier~\cite{felder}, we shall construct outgoing
intertwining vectors in the same way as
\cite{hasegawa1,hasegawa2,quano}.
We can consequently define factorized
L-operators $\Lcheckk(\x)$ (Theorem~\ref{thm:loperator}).
\begin{thm}[Factorized L-operator]
For $\x_1,\x_2\not\in\Bbb Z + \Bbb Z \t$
\begin{equation*}
(1\otimes \Rcheck(\x_{12}))(\Lcheckk(\x_{1})\otimes
1)(1\otimes \Lcheckk(\x_{2}))
=(\Lcheckk(\x_{2})\otimes 1)(1\otimes
\Lcheckk(\x_{1}))(\Rcheck(\x_{12})\otimes 1).
  \end{equation*}
\end{thm}
In the last section, after stating the results obtained
by Felder and Pasquier \cite{felder} more precisely,
we show that the vertex--IRF correspondence and the factorized
L-operators for the
elliptic $R$-operator
imply those for Belavin's
$R$-matrix.
\section{Review of the properties of an elliptic $R$-operator}
In this section, we review the construction and the
properties of an elliptic
$R$-operator
\cite{felder,shibukawa1,shibukawa2,shibukawa3}.
We fix $\t\in \Bbb C$ such that $\mathrm{Im}\ \t >0$ and
define an open subset $D \subset \Bbb C$ by
\[ D=\{ z\in \Bbb C \ ; \ |\mathrm{Im}\ z|<\frac{\mathrm{Im}\ \t}2 \}.
\]
Let $\V$ be a space of all functions $f$ holomorphic on $D$ and
such that
\[ f(z+1)=f(z) \quad \forall z\in D.
\]
Similarly let $\V \otimeshat \V$ be a space of all functions
$f$ holomorphic on $D\times D$ with the property
\[ f(z_1+1,z_2)=f(z_1,z_2+1)=f(z_1,z_2) \quad
\forall z_1,z_2\in D.
\]

Now we define an elliptic $R$-operator $\Rcheck(\x)$ on $\V
\otimeshat \V$.
Let $\m$ be a complex number such that $\m \not\in \Bbb
Z + \Bbb Z \t$ and let $\vt1(z)=\vt1(z,\t)$ be an
elliptic theta
function
\[ \vt1(z)=\sum_{m\in \Bbb Z}\exp [\p \rt-1 (m+\frac12)^2\t
+2\p \rt-1 (m+\frac12)(z+\frac12)].
\]
The elliptic theta function $\vt1(z)$ satisfies the following
properties.
\begin{enumerate}
\item $\vt1(z)$ is entire,
\item \(\vt1(z+1)=-\vt1(z),\)
\item \(\vt1(z+\t)=-\exp (-2\p\rt-1 z-\p\rt-1\t) \vt1(z),\)
\item $\vt1(z)$ has simple zeros at $z\in \Bbb Z + \Bbb Z
  \t$,
\item $\vt1(z)$ satisfies the three term equation
(cf.~\cite{whittaker} p.461);
\begin{align*}
&\phantom{+\vartheta}\vt1(x+y)\vt1(x-y)\vt1(z+w)\vt1(z-w)\\
&+\vt1(x+z)\vt1(x-z)\vt1(w+y)\vt1(w-y)\\
&+\vt1(x+w)\vt1(x-w)\vt1(y+z)\vt1(y-z)\\
&=0,
\end{align*}
\item $\vt1(-z)=-\vt1(z).$
\end{enumerate}
\begin{defn}[Elliptic ${\bold R}$-operator]\label{defn:roperator}
For $f\in \V\otimeshat \V$, we define
\begin{equation*}
(\Rcheck (\x )f)(z_1,z_2):=
\frac{\vt1(\x)\vt1(z_{21}-\m)\vt1'(0)}
{\vt1(-\m)\vt1(z_{21})}f(z_2,z_1)
+\frac{\vt1(z_{21}-\x)\vt1'(0)}{\vt1(z_{21})}
f(z_1,z_2),
\end{equation*}
where $z_{21}:=z_2-z_1$,
$\vt1'(0)=\frac{d\vt1}{dz}(0,\t)$
 and $\x \in \Bbb C$.
The complex number $\x$ is called a spectral parameter.
\end{defn}
We set $X=\{(z_1,z_2)\in D\times D \ ; \ z_{21}\in \Bbb Z \}$.
By the property (4) of the elliptic theta function $\vt1(z)$,
the function $\Rcheck(\x)f$ has the singularities at the
points $(z_1,z_2) \in X$.
The lemma below tells us that all singularities are
removable.
\begin{lem}
  There is a unique function $F$ holomorphic on $D\times D$
  and such that $F(z_1,z_2)=(\Rcheck (\x )f)(z_1,z_2)$ for
  $(z_1,z_2)  \in D \times D \backslash X$.
\end{lem}
\begin{pf}
For $(z_1,z_2)\in D\times D \backslash X$ and $m\in \Bbb Z$,
\begin{align*}
  &(\Rcheck (\x )f)(z_1,z_2)\\
=&
\frac{\vt1(\x)\vt1'(0)f(z_2,z_1)}{\vt1(-\m)}\cdot
\frac{\vt1(z_{21}-\m-m)-\vt1(-\m)}{z_2-z_1-m}\cdot
\frac{z_2-z_1-m}{\vt1(z_{21}-m)}\\
+&
\vt1(\x)\vt1'(0)\frac{f(z_2-m,z_1)-f(z_1,z_1)}{z_2-z_1-m}\cdot
\frac{z_2-z_1-m}{\vt1(z_{21}-m)}\\
+&
\vt1'(0)
\frac{z_2-z_1-m}{\vt1(z_{21}-m)}
\{
f(z_1,z_1)
\frac{\vt1(z_{21}-\x-m)-\vt1(-\x)}{z_2-z_1-m}\\
&+
\vt1(z_{21}-\x-m)
\frac{f(z_1,z_2-m)-f(z_1,z_1)}{z_2-z_1-m}\}.
\end{align*}
Thus there is a function $F$ continuous
on $D\times D$ and such that $F(z_1,z_2)=(\Rcheck (\x
)f)(z_1,z_2)$ for $(z_1,z_2)\in D\times D \backslash X$.
In fact, we define
\begin{align*}
  &F(z_1,z_2)\\
=
&
\left\{
\begin{array}{ll}
\frac{\vt1(\x)\vt1'(-\m)+\vt1'(-\x)\vt1(-\m)}{\vt1(-\m)}f(z,z)
+\vt1(\x)
(
\frac{\partial f}{\partial z_1}(z,z)
-\frac{\partial f}{\partial z_2}(z,z)
),
&
(z_1,z_2)=(z,z+m), \\
(\Rcheck (\x )f)(z_1,z_2),
&
\text{otherwise.}
\end{array}
\right.
\end{align*}
Making use of the
Riemann removable singularity theorem (cf.~\cite{gunning}),
this function $F$ is holomorphic on $D\times D$.
\end{pf}
We also denote by $\Rcheck (\x)f$ this holomorphic function
$F$.
It is easy to see that
\[
(\Rcheck (\x )f)(z_1+1,z_2)=(\Rcheck (\x )f)(z_1,z_2+1)=
(\Rcheck (\x )f)(z_1,z_2)
\]
for $(z_1,z_2) \in D\times D$.
Hence
\(
\Rcheck(\x)f \in \V \otimeshat \V
\)
for $f \in \V \otimeshat \V$,
and $\Rcheck(\x)$ is an operator on $\V \otimeshat
\V$ as a result.

Let $\V\otimeshat \V\otimeshat \V$ be a space of all
functions $f$ on $D\times D\times D$ and such that
\[
f(z_1+1,z_2,z_3)=f(z_1,z_2+1,z_3)=f(z_1,z_2,z_3+1)=f(z_1,z_2,z_3)
\quad \forall z_1,z_2,z_3\in D.
\]
By the three term equation of $\vt1(z)$ (the property (5)),
we get the following theorem.
\begin{thm}[\cite{shibukawa1,shibukawa2,shibukawa3}]
\label{thm:shibukawa}
$\Rcheck(\x)$ satisfies the Yang--Baxter equation on
$\V\otimeshat \V\otimeshat \V$\rom;
\begin{equation}
(1\otimes \Rcheck(\x_{12}))(\Rcheck(\x_{13})\otimes
1)(1\otimes \Rcheck(\x_{23}))
=(\Rcheck(\x_{23})\otimes 1)(1\otimes
\Rcheck(\x_{13}))(\Rcheck(\x_{12})\otimes 1),
\label{intertwine:ybe}
\end{equation}
where $\x_{ij}=\x_i-\x_j$.
\end{thm}
For $\x \in \Bbb C$ and $n=1,2,\ldots$, let $V_n(\x)$ be a
space of all functions $f$ holomorphic on $\Bbb C$ and such
that
\begin{align*}
  f(z+1)&=f(z),\\
  f(z+\t)&=(-1)^n\exp (2\p\rt-1 (\x-nz))f(z).
\end{align*}
It is well known that $V_n(\x)$ has dimension $n$.
We easily see that
\begin{equation}
\{ \vartheta {\frac12-\frac{j}{n} \brack{\frac{n}{2}}}(\x-nz,n\t)
\exp (\p\rt-1 nz) \}_{j\in \Bbb Z/n\Bbb Z}
\label{summary:basis}
\end{equation}
is a basis of $V_n(\x)$.
Here
\(
\vartheta {a \brack{b}}(z,\t)
\)
is a theta function with rational characteristics;
\[
\vartheta {a \brack{b}}(z,\t)=
\sum_{m\in\Bbb Z}\exp [\p\rt-1 (m+a)^2\t+2\p\rt-1 (m+a)(z+b)].
\]
In \cite{felder} Felder and Pasquier show the following.
\begin{thm}[\cite{felder}]\label{thm:felder}
\(
\Rcheck(\x_{12})(V_n(\x_1)\otimes V_n(\x_2+\m))
\subset V_n(\x_2)\otimes V_n(\x_1+\m).
\)
\end{thm}
\begin{rem}\label{rem:roperatorvm}
Let $\Vm$ be a space of all functions $f$ holomorphic on $D$ and such
that
\[
f(z+1)=-f(z).
\]
We set $\Vm\otimeshat\Vm$ and
$\Vm\otimeshat\Vm\otimeshat\Vm$ in the same way as $\V$.
Then we can define the elliptic $R$-operator $\Rcheck(\x)$ on
$\Vm\otimeshat\Vm$, which is the same as in
Definition~\ref{defn:roperator}.
It is easy to see that $\Rcheck(\x)$ on
$\Vm\otimeshat\Vm$ satisfies the Yang--Baxter equation
\eqref{intertwine:ybe}.

We denote $V_n^-(\x)$ as a space of all functions $f$
holomorphic on $\Bbb C$ and such that
\begin{align*}
  f(z+1)&=-f(z),\\
  f(z+\t)&=(-1)^n\exp2\p\rt-1 (\x-nz+\frac{\t}{2}) f(z).
\end{align*}
We have
\[
\Rcheck(\x_{12})(V_n^-(\x_1)\otimes V_n^-(\x_2+\m))
\subset V_n^-(\x_2)\otimes V_n^-(\x_1+\m).
\]
A basis of $V_n^-(\x)$ is as follows.
\[
\{ \vartheta {\frac12-\frac{j}{n} \brack{\frac{n}{2}}}(\x-nz,n\t)
\exp (\p\rt-1 (n+1)z) \}_{j\in \Bbb Z/n\Bbb Z}.
\]
\end{rem}
\begin{rem}\label{rem:roperator}
  Let $\cal M$ be a space of the meromorphic functions on
  $\Bbb C^2$. Then we note that the elliptic $R$-operator
  $\Rcheck(\x)$ can be regarded as an operator on $\cal M$
  and satisfies the Yang--Baxter
  equation \eqref{intertwine:ybe}.
\end{rem}
\section{Incoming intertwining vectors and vertex--irf
  correspondence}
In what follows $\m \in \Bbb R \backslash \Bbb Z$,
and let $\L$ be a set of sequences $\l=(\l_i)\;\;(i\in\Bbb
  Z)$ such that
\begin{align*}
&\l_i\in D,\\
&\l_{ij}:=\l_i-\l_j
\not\in \Bbb Z +\Bbb Z \m \quad \forall i\ne j\in \Bbb Z.
\end{align*}
We take $r \in \Bbb R$ such that
\(
r \not\in \Bbb Q +\Bbb Q \m,
\)
and set
\[
\h_i:=ir \quad (i\in \Bbb Z).
\]
Then $\h=(\h_i) \in \L$. Hence, for any $\m$, the set
$\L$ is not empty.
For $i\in \Bbb Z$, we define the sequences
$\ve_i=(\d_{ij})\;\;(j\in \Bbb Z)$,
and for $\l \in \L$, let $\l+\m\ve_i$ denote the sequence
\begin{equation*}
  (\l+\m\ve_i)_j=
  \begin{cases}
    \l_j, & j \neq i,\\
    \l_i+\m, & j=i.
  \end{cases}
\end{equation*}
We note that
$\l+\m\ve_i\in \L$ for all $i\in \Bbb Z$
by the definition of $\L$.
\begin{defn}[Boltzmann weight of the IRF model]
\label{defn:boltzmann}
For $\l,\k,\k',\n \in \L$,
Boltzmann weights
\(
\Wcheck
\begin{bmatrix}
  &\k'&\\
\l&\x&\n\\
&\k&
\end{bmatrix}
\in \Bbb C
\)
of an interaction--round--a--face (IRF)
model are given as follows
(cf.~\cite{baxter,hasegawa1,hasegawa2,jimbo,quano}).
For $\l\in\L$, we put
\begin{align*}
\Wcheck
\begin{bmatrix}
  &\l+\m\ve_i&\\
\l&\x&\l+2\m\ve_i\\
&\l+\m\ve_i&
\end{bmatrix}
&:=\frac{\vt1(\m-\x)\vt1'(0)}{\vt1(\m)},\\
\Wcheck
\begin{bmatrix}
  &\l+\m\ve_i&\\
\l&\x&\l+\m(\ve_i+\ve_j)\\
&\l+\m\ve_i&
\end{bmatrix}
&:=\frac{\vt1(\l_{ji}-\x)\vt1'(0)}{\vt1(\l_{ji})} \quad (i\ne
j),\\
\Wcheck
\begin{bmatrix}
  &\l+\m\ve_j&\\
\l&\x&\l+\m(\ve_i+\ve_j)\\
&\l+\m\ve_i&
\end{bmatrix}
&:=\frac{\vt1(\x)\vt1(\l_{ji}-\m)\vt1'(0)}{\vt1(\l_{ji})\vt1(-\m)}
\quad (i\ne j),
\end{align*}
otherwise we set
\[
\Wcheck
\begin{bmatrix}
  &\k'&\\
\l&\x&\n\\
&\k&
\end{bmatrix}
:=0.
\]
\end{defn}
Next we define incoming intertwining vectors
of the elliptic
$R$-operator.
\begin{defn}[Incoming intertwining vector]
\label{defn:incoming}
For $\l,\k\in\L$,
define an incoming intertwining vector
$\pbar_\l^\k\in \cal
V^*$
as follows;
\begin{equation*}
  \pbar_\l^\k f:=
  \begin{cases}
    f(\l_i),&\exists i\in\Bbb Z \;\mathrm{s.t.}\;\k=\l+\m\ve_i,\\
    0,&\text{otherwise.}
  \end{cases}
\end{equation*}
\end{defn}
The incoming intertwining vectors are the Dirac delta
functions essentially.
By Definition~\ref{defn:roperator}
we can get a vertex--IRF correspondence for the elliptic
$R$-operator.
\begin{thm}[Vertex--IRF correspondence]\label{thm:vertex-irf}
For $\l,\k,\n\in\L$
\begin{equation}
\pbar_\l^\k \otimes \pbar_\k^\n \Rcheck(\x)
=\sum_{\k'\in\L}
\Wcheck
\begin{bmatrix}
  &\k'&\\
\l&\x&\n\\
&\k&
\end{bmatrix}
\pbar_\l^{\k'}\otimes \pbar_{\k'}^\n,
\label{vertex-irf:vertex-irf}
\end{equation}
where both sides are the operators
\(
\V \otimeshat \V \to \Bbb C
\).
\end{thm}
It is to be noted that, by Definition~\ref{defn:boltzmann}
and \ref{defn:incoming},
the both sides of equation
\eqref{vertex-irf:vertex-irf} are zero unless there
exist $i,j \in \Bbb Z$ such that $\k=\l+\m\ve_i,\
\n=\l+\m(\ve_i+\ve_j)$.
The other cases are as follows.
\begin{align*}
 \pbar_\l^{\l+\m\ve_i} \otimes
 \pbar_{\l+\m\ve_i}^{\l+2\m\ve_i} \Rcheck(\x)
&=\frac{\vt1(\m-\x)\vt1'(0)}{\vt1(\m)}
 \pbar_\l^{\l+\m\ve_i} \otimes
 \pbar_{\l+\m\ve_i}^{\l+2\m\ve_i},\\
\\
 \pbar_\l^{\l+\m\ve_i} \otimes
 \pbar_{\l+\m\ve_i}^{\l+\m(\ve_i+\ve_j)} \Rcheck(\x)
&=\frac{\vt1(\l_{ji}-\x)\vt1'(0)}{\vt1(\l_{ji})}
 \pbar_\l^{\l+\m\ve_i} \otimes
 \pbar_{\l+\m\ve_i}^{\l+\m(\ve_i+\ve_j)} \\
&+\frac{\vt1(\x)\vt1(\l_{ji}-\m)\vt1'(0)}{\vt1(\l_{ji})\vt1(-\m)}
 \pbar_\l^{\l+\m\ve_j} \otimes
 \pbar_{\l+\m\ve_j}^{\l+\m(\ve_i+\ve_j)},
\end{align*}
for $i \neq j$.
Since $\Rcheck(\x)$ satisfies the Yang--Baxter equation
\eqref{intertwine:ybe},
we can show
\begin{cor}\label{cor:star-triangle}
  The Boltzmann weights of the IRF model satisfy the
  star--triangle relation\rom;
  \begin{equation}
    \begin{split}
      \sum_{\k'\in\L}
&
\Wcheck
\begin{bmatrix}
  &\k'&\\
\k&\x_{12}&\g\\
&\n&
\end{bmatrix}
\Wcheck
\begin{bmatrix}
  &\a&\\
\l&\x_{13}&\k'\\
&\k&
\end{bmatrix}
\Wcheck
\begin{bmatrix}
  &\b&\\
\a&\x_{23}&\g\\
&\k'&
\end{bmatrix}
\\
& \qquad =
\sum_{\k'\in\L}
\Wcheck
\begin{bmatrix}
  &\k'&\\
\l&\x_{23}&\n\\
&\k&
\end{bmatrix}
\Wcheck
\begin{bmatrix}
  &\b&\\
\k'&\x_{13}&\g\\
&\n&
\end{bmatrix}
\Wcheck
\begin{bmatrix}
  &\a&\\
\l&\x_{12}&\b\\
&\k'&
\end{bmatrix}
,
    \end{split}
\label{vertex-irf:star-triangle}
  \end{equation}
for $\l,\k,\n,\a,\b,\g\in\L$.
\end{cor}
\begin{pf}
Unless there exist $i,j,k \in \Bbb Z$ such that
$\k=\l+\m\ve_i, \ \n=\l+\m(\ve_i+\ve_j)$ and
$\g=\l+\m(\ve_i+\ve_j+\ve_k)$, the both sides of
equation \eqref{vertex-irf:star-triangle} are zero.
Then we assume that
\[
\k=\l+\m\ve_i, \ \n=\l+\m(\ve_i+\ve_j), \
\g=\l+\m(\ve_i+\ve_j+\ve_k)\quad (i,j,k \in \Bbb Z).
\]
Moreover the both sides of equation
\eqref{vertex-irf:star-triangle} are
zero unless
\begin{align*}
&\a=\l+\m\ve_i, \  \l+\m\ve_j \ \text{or} \ \l+\m\ve_k\\
\intertext{and}
&\b=\l+\m(\ve_i+\ve_j), \ \l+\m(\ve_i+\ve_k) \ \text{or}\
\l+\m(\ve_j+\ve_k),
\end{align*}
so it suffices to show equation
\eqref{vertex-irf:star-triangle} in each case.

Since $\Rcheck(\x)$ satisfies the Yang--Baxter equation
\eqref{intertwine:ybe},
\begin{equation*}
\begin{split}
&((1\otimes \Rcheck(\x_{12}))(\Rcheck(\x_{13})\otimes
1)(1\otimes \Rcheck(\x_{23}))f)(z_1,z_2,z_3)\\
& \qquad
=((\Rcheck(\x_{23})\otimes 1)(1\otimes
\Rcheck(\x_{13}))(\Rcheck(\x_{12})\otimes 1)f)(z_1,z_2,z_3).
\end{split}
\end{equation*}
Putting $z_1=\l_i, \ z_2=\l_j+\m\d_{ij}$ and
$z_3=\l+\m(\d_{ik}+\d_{jk})$ in the coefficient of
$f(z_1,z_2,z_3)$, we obtain equation
\eqref{vertex-irf:star-triangle} in the case $\a=\l+\m\ve_i$
and $\b=\l+\m(\ve_i+\ve_j)$.
We can prove the other cases in the similar way, so we omit
the proof.
\end{pf}
\begin{rem}\label{rem:incomingvm}
We define an incoming intertwining vector
$\pbar_\l^\k\in (\Vm)^*$ in the same way as
Definition~\ref{defn:incoming}; for $f\in\Vm$,
\begin{equation*}
  \pbar_\l^\k f:=
  \begin{cases}
    f(\l_i),&\exists i\in\Bbb Z \;\mathrm{s.t.}\;\k=\l+\m\ve_i,\\
    0,&\text{otherwise.}
  \end{cases}
\end{equation*}
In this case, we also get a vertex--IRF correspondence;
for $\l,\k,\n\in\L$
\[ \pbar_\l^\k \otimes \pbar_\k^\n \Rcheck(\x)
=\sum_{\k'\in\L}
\Wcheck
\begin{bmatrix}
  &\k'&\\
\l&\x&\n\\
&\k&
\end{bmatrix}
\pbar_\l^{\k'}\otimes \pbar_{\k'}^\n.
\]
\end{rem}
\section{Outgoing intertwining vectors and factorized L-operators}
To begin with, we define outgoing intertwining vectors of the
elliptic $R$-operator (cf. \cite{hasegawa1,hasegawa2,quano}).

Let $k_1$ and $k_2$ be integers such that $k_1 \leq k_2$, and we
set $\bk :=(k_1,k_2)$ and $k=k_2-k_1+1$.
For $\l,\k\in\L$ and $k_1\leq j\leq k_2$,
we define $\pbar_{\bk}(\x)_\l^{\k\, j}\in\Bbb C$ by
\[
  \pbar_{\bk}(\x)_\l^{\k\,j}:=
    \pbar_\l^\k(\vartheta {\frac12-\frac{j-k_1}k \brack{\frac{k}2}}
(\x+|\l|_{\bk}-kz,k\t)\exp (\p\rt-1 kz)),
\]
where
$|\l|_{\bk}=\sum_{i=k_1}^{k_2}\l_i$.
\begin{prop}\label{prop:weyl-kac}
  For $\l\in\L$ and $\x\not\in \Bbb
Z + \Bbb Z \t$,
the $k$--by--$k$ matrix
\(
\begin{pmatrix}
  \pkbar(\x)_{\l}^{\l+\m\ve_i\, j}
\end{pmatrix}_{k_1\leq i,j\leq k_2}
\)
is invertible.
\end{prop}
\begin{pf}
Since
\begin{align*}
&\begin{pmatrix}
  \pkbar(\x)_{\l}^{\l+\m\ve_i,\ j}
\end{pmatrix}_{k_1\leq i,j\leq k_2}\\
=&
\mathrm{diag}(\exp\p\rt-1 k\l_{k_1},\ldots,\exp\p\rt-1
k\l_{k_2})
\begin{pmatrix}
  \vartheta {\frac12-\frac{j-k_1}k \brack{\frac{k}2}}
(\x+|\l|_{\bk}-k\l_i, k\t)
\end{pmatrix}_{k_1\leq i,j\leq k_2},
\end{align*}
it suffices to prove
\[
\det
\begin{pmatrix}
  \vartheta {\frac12-\frac{j-k_1}k \brack{\frac{k}2}}
(\x+|\l|_{\bk}-k\l_i, k\t)
\end{pmatrix}_{k_1\leq i,j\leq k_2}
\ne 0.
\]
The Weyl--Kac denominator formula for $A_{k-1}^{(1)}$
(cf.~\cite{kac,hasegawa2}) yields
\begin{align*}
&\det
\begin{pmatrix}
  \vartheta {\frac12-\frac{j}k \brack{\frac{k}2}}
(ku_i, k\t)
\end{pmatrix}_{1\leq i,j\leq k}\\
=&
(\rt-1 \h(\t))^{-\frac12(k-1)(k-2)}\vt1(\sum_{i=1}^ku_i)
\prod_{1\leq j<i\leq k}\vt1(u_{ij}).
\end{align*}
Here $\h(\t)$ is Dedekind's $\h$--function
\[
\h(\t)=\exp\frac{\p\rt-1 \t}{12}\prod_{m=1}^\infty
(1-\exp2\p\rt-1 m\t).
\]
Then we obtain
\begin{align*}
&  \det
\begin{pmatrix}
  \vartheta {\frac12-\frac{j-k_1}k \brack{\frac{k}2}}
(\x+|\l|_{\bk}-k\l_i, k\t)
\end{pmatrix}_{k_1\leq i,j\leq k_2}
\\
=&
(-1)^{k-1}
\begin{pmatrix}
  \vartheta {\frac12-\frac{j}k \brack{\frac{k}2}}
(\x+|\l|_{\bk}-k\l_{i+k_1-1}, k\t)
\end{pmatrix}_{1\leq i,j\leq k}\\
=&
(-1)^{k-1}(\rt-1 \h(\t))^{-\frac12 (k-1)(k-2)}
\vt1(\x)\prod_{k_1\leq i<j\leq k_2}\vt1(\l_{ij}),
\end{align*}
thereby completing the proof.
\end{pf}
The proposition above says that for $\l, \k\in\L$, $k_1 \leq
j \leq k_2$ and
$\x\not\in \Bbb Z + \Bbb Z \t$,
there exist
$\pk(\x)_{\l\, j}^\k \in\Bbb C$ which are characterized by the
  following duality relations;
\begin{equation}
\left\{
  \begin{aligned}
    &\sum_{i=k_1}^{k_2}\pk(\x)_{\l\,j}^{\l+\m\ve_i}
\pkbar(\x)_\l^{\l+\m\ve_i\,\ell}
=\d_{j\ell},\\
&\sum_{i=k_1}^{k_2}\pkbar(\x)_{\l}^{\l+\m\ve_j\,i}
\pk(\x)_{\l\,i}^{\l+\m\ve_\ell}
=\d_{j\ell},
  \end{aligned}
\right.\label{outgoing:dualityrelation}
\end{equation}
and for $\k\ne\l+\m\ve_i\;(k_1\leq \forall i\leq k_2)$ we set
\[
  \pk(\x)_{\l\ j}^\k:=0.
\]
\begin{defn}[Outgoing intertwining vector]\label{defn:outgoing}
For $\l,\k \in \L$ and $\x \not\in \Bbb Z + \Bbb Z \t$, an
outgoing intertwining vector
\(
\pk(\x)_\l^\k(z)\in V_k(\x+|\l|_{\bk})
\)
of the elliptic $R$-operator is defined as follows
(cf.~\eqref{summary:basis}).
\[
  \pk(\x)_\l^\k(z):=
    \sum_{j=k_1}^{k_2}\pk(\x)_{\l\, j}^\k
\vartheta {\frac12-\frac{j-k_1}k \brack{\frac{k}2}}
(\x+|\l|_{\bk}-kz,k\t)\exp (\p\rt-1 kz).
\]
\end{defn}
Equation \eqref{outgoing:dualityrelation} is equivalent to
\begin{equation}
  \left\{
\begin{aligned}
  &\sum_{i=k_1}^{k_2}\pk(\x)_\l^{\l+\m\ve_i}(z)
\pbar_\l^{\l+\m\ve_i}=\mathrm{id}\quad \text{on}\
V_k(\x+|\l|_{\bk}),\\
&\pbar_\l^{\l+\m\ve_i}(\pk(\x)_\l^{\l+\m\ve_j})
=\d_{ij} \quad \text{for} \ k_1 \leq i,j\leq k_2.
\end{aligned}
\right. \label{outgoing:dualityrelation2}
\end{equation}
The outgoing intertwining vectors satisfy the following.
\begin{prop}\label{prop:outgoing}
For $\l,\k,\n\in\L$ and $\x_1,\x_2\not\in\Bbb
Z + \Bbb Z \t$,
  \[
(\Rcheck(\x_{12})\pk(\x_1)_\l^\k \otimes
\pk(\x_2)_\k^\n)(z_1,z_2)=
\sum_{\k'\in\L}
\pk(\x_2)_\l^{\k'}(z_1)\otimes
\pk(\x_1)_{\k'}^\n(z_2)
\Wcheck
\begin{bmatrix}
  &\k&\\
\l&\x_{12}&\n\\
&\k'&
\end{bmatrix}.
\]
\end{prop}
\begin{pf}
By Definition~\ref{defn:boltzmann} and \ref{defn:outgoing},
it suffices to show
\begin{align*}
&(\Rcheck(\x_{12})\pk(\x_1)_\l^{\l+\m\ve_i} \otimes
\pk(\x_2)_{\l+\m\ve_i}^{\l+\m(\ve_i+\ve_j)})(z_1,z_2)\\
=&
\sum_{\ell =k_1}^{k_2}
\pk(\x_2)_\l^{\l+\m\ve_{\ell}}(z_1)\otimes
\pk(\x_1)_{\l+\ve_{\ell}}^{\l+\m(\ve_i+\ve_j)}(z_2)
\Wcheck
\begin{bmatrix}
  &\l+\m\ve_i&\\
\l&\x_{12}&\l+\m(\ve_i+\ve_j)\\
&\l+\m\ve_{\ell}&
\end{bmatrix}
\end{align*}
for any $\l\in\L$ and $k_1\leq \forall i,j\leq k_2$.
With the aid of Theorem~\ref{thm:vertex-irf} and
equation~\eqref{outgoing:dualityrelation2}, we obtain for
$k_1 \leq \forall a,b\leq k_2$,
\begin{align*}
&
\pbar_\l^{\l+\m\ve_i} \otimes
\pbar_{\l+\m\ve_i}^{\l+\m(\ve_i+\ve_j)}
((\Rcheck(\x_{12})
\pk(\x_1)_\l^{\l+\m\ve_a}\otimes
\pk(\x_2)_{\l+\m\ve_a}^{\l+\m(\ve_a+\ve_b)})(z_1,z_2))\\
=&
\sum_{\ell=k_1}^{k_2}
\Wcheck
\begin{bmatrix}
  &\l+\m\ve_{\ell}&\\
\l&\x_{12}&\l+\m(\ve_i+\ve_j)\\
&\l+\m\ve_i&
\end{bmatrix}\\
\phantom{=}&\times
(\pbar_\l^{\l+\m\ve_{\ell}}\otimes
\pbar_{\l+\m\ve_{\ell}}^{\l+\m(\ve_i+\ve_j)})
(\pk(\x_1)_\l^{\l+\m\ve_a}(z_1)\otimes
\pk(\x_2)_{\l+\m\ve_a}^{\l+\m(\ve_a+\ve_b)}(z_2))\\
=&
\Wcheck
\begin{bmatrix}
  &\l+\m\ve_a&\\
\l&\x_{12}&\l+\m(\ve_i+\ve_j)\\
&\l+\m\ve_i&
\end{bmatrix}
\d_{\l+\m(\ve_i+\ve_j)\ \l+\m(\ve_a+\ve_b)}.
\end{align*}
Then
\begin{align*}
 &\sum_{i,j=k_1}^{k_2}(\pk(\x_2)_\l^{\l+\m\ve_i}(z_1)
\otimes \pk(\x_1)_{\l+\m\ve_i}^{\l+\m(\ve_i+\ve_j)}(z_2))
\\
\phantom{=}&\times
(\pbar_\l^{\l+\m\ve_i} \otimes
\pbar_{\l+\m\ve_i}^{\l+\m(\ve_i+\ve_j)})
((\Rcheck(\x_{12})
\pk(\x_1)_\l^{\l+\m\ve_a} \otimes
\pk(\x_2)_{\l+\m\ve_a}^{\l+\m(\ve_a+\ve_b)})(z_1,z_2))\\
=&
 \sum_{i,j=k_1}^{k_2}\pk(\x_2)_\l^{\l+\m\ve_i}(z_1)
\otimes \pk(\x_1)_{\l+\m\ve_i}^{\l+\m(\ve_i+\ve_j)}(z_2)
\\
\phantom{=}&\times
\Wcheck
\begin{bmatrix}
  &\l+\m\ve_a&\\
\l&\x_{12}&\l+\m(\ve_i+\ve_j)\\
&\l+\m\ve_i&
\end{bmatrix}
\d_{\l+\m(\ve_i+\ve_j)\ \l+\m(\ve_a+\ve_b)}\\
=&
\sum_{i=k_1}^{k_2}\pk(\x_2)_\l^{\l+\m\ve_i}(z_1)
\otimes \pk(\x_2)_{\l+\m\ve_i}^{\l+\m(\ve_a+\ve_b)}(z_2)
\Wcheck
\begin{bmatrix}
  &\l+\m\ve_a&\\
\l&\x_{12}&\l+\m(\ve_a+\ve_b)\\
&\l+\m\ve_i&
\end{bmatrix}.
\end{align*}
By virtue of Definition~\ref{defn:outgoing} and
Theorem~\ref{thm:felder} we deduce
\[
(\Rcheck(\x_{12})
\pk(\x_1)_\l^{\l+\m\ve_a} \otimes
\pk(\x_2)_{\l+\m\ve_a}^{\l+\m(\ve_a+\ve_b)})(z_1,z_2)
\in V_k(\x_2+|\l|_{\bk})\otimes V_k(\x_1+|\l|_{\bk}+\m).
\]
{}From equation~\eqref{outgoing:dualityrelation2}, we are led
to the desired result.
\end{pf}

For $\l,\k \in \L$ and $\x\not\in\Bbb Z + \Bbb Z \t$, we define
an operator
\(\Lcheckk(\x)_\l^\k : \V\to\V \)
by
\[ (\Lcheckk(\x)_\l^\k f)(z):= \pk(\x)_\l^\k(z)\pbar_\l^\k f
\quad (f \in \V).
\]
Theorem~\ref{thm:vertex-irf} and Proposition~\ref{prop:outgoing}
say
\begin{lem}\label{lem:loperator}
For $\l,\n\in\L$ and
$\x_1,\x_2\not\in\Bbb Z + \Bbb Z \t$,
  \begin{equation*}
    \sum_{\k\in\L}\Rcheck(\x_{12})\Lcheckk(\x_1)_\l^\k
\otimes \Lcheckk(\x_2)_\k^\n
=\sum_{\k\in\L}
\Lcheckk(\x_2)_\l^\k \otimes \Lcheckk(\x_1)_\k^\n
\Rcheck(\x_{12})
  \end{equation*}
on $\V\otimeshat\V$.
\end{lem}
\begin{pf}
For $f \in \V \otimeshat \V$,
\begin{align*}
&\sum_{\k\in\L}(\Rcheck(\x_{12})\Lcheckk(\x_1)_\l^\k
\otimes \Lcheckk(\x_2)_\k^\n f)(z_1,z_2)\\
=&
\sum_{\k\in\L}(\Rcheck(\x_{12})\pk(\x_1)_\l^\k
\otimes \pk(\x_2)_\k^\n)(z_1,z_2) \cdot (\pbar_\l^\k \otimes
\pbar_\k^\n)f\\
=&
\sum_{\k,\k'\in\L}\pk(\x_2)_\l^{\k'}(z_1)
\otimes \pk(\x_1)_{\k'}^\n(z_2)
\Wcheck
\begin{bmatrix}
  &\k&\\
\l&\x_{12}&\n\\
&\k'&
\end{bmatrix}
(\pbar_\l^\k \otimes
\pbar_\k^\n)f\\
=&
\sum_{\k'\in\L}\pk(\x_2)_\l^{\k'}(z_1)
\otimes \pk(\x_1)_{\k'}^\n(z_2)
(\pbar_\l^{\k'} \otimes
\pbar_{\k'}^\n \Rcheck(\x_{12})f)\\
=&
\sum_{\k\in\L}
(\Lcheckk(\x_2)_\l^\k \otimes \Lcheckk(\x_1)_\k^\n
\Rcheck(\x_{12}) f)(z_1,z_2).
\end{align*}
We have thus proved the lemma.
\end{pf}

Now we are in the position to construct factorized
L-operators for the elliptic $R$-operator.
Let $\W$ be a space of all $\Bbb C$--valued
functions on $\L$,
and let
$\V\otimeshat \W$ (resp.~$\W \otimeshat \V$)
be a
space of all functions $g : D\times \L\to \Bbb C$
(resp.~$\L\times D \to \Bbb C$)
such that
$g(\cdot,\l)\in\V$ (resp.~$g(\l,\cdot)\in\V$)
for any $\l\in\L$.
We define a factorized L-operator
\(
\Lcheckk(\x):\V \otimeshat \W \to
\W \otimeshat \V \) as follows
\cite{bazhanov,hasegawa1,hasegawa2,quano}.
For $g\in \V\otimeshat \W$ and $\x\not\in \Bbb
Z + \Bbb Z \t$
\begin{equation}
  (\Lcheckk(\x)g)(\l,z):=
\sum_{\k\in\L}(\Lcheckk(\x)_\l^\k g(\cdot,\k))(z).
\label{loperator:loperator}
\end{equation}
For $\l\in\L$ we set $\d^\l\in\W$ as follows;
\[
\d^\l(\k)=\d_{\l\k}.
\]
We note that \(\W=\prod_{\k\in\L}\Bbb C\d^\k\)
(cf.~\cite{hasegawa1}).
Then, for $f\in\V$,
\begin{equation*}
( \Lcheckk ( \x )( f \otimes \d^\k ))(\l,z)
=(\Lcheckk(\x)_\l^\k f)(z),
\end{equation*}
and equation~\eqref{loperator:loperator} is hence equivalent to
\begin{equation*}
\Lcheckk(\x)(f\otimes \d^\k)
= \sum_{\l\in\L}\d^\l\otimes \Lcheckk(\x)_\l^\k f.
\end{equation*}

We define $\V\otimeshat\V\otimeshat \W$
(resp.~$\W \otimeshat\V\otimeshat\V$) by a
space of all functions $g : D\times D\times\L\to \Bbb C$
(resp.~$\L\times D\times D \to \Bbb C$)
such that
$g(\cdot,\cdot,\l)\in\V\otimeshat \V$
(resp.~$g(\l,\cdot,\cdot)\in\V\otimeshat \V$)
for any $\l\in\L$.
By means of Lemma \ref{lem:loperator}, we immediately obtain
the following theorem.
\begin{thm}[Factorized L-operator]\label{thm:loperator}
For $\x_1,\x_2\not\in\Bbb Z + \Bbb Z \t$
  \begin{equation*}
(1\otimes \Rcheck(\x_{12}))(\Lcheckk(\x_{1})\otimes
1)(1\otimes \Lcheckk(\x_{2}))
=(\Lcheckk(\x_{2})\otimes 1)(1\otimes
\Lcheckk(\x_{1}))(\Rcheck(\x_{12})\otimes 1),
  \end{equation*}
where both sides are the operators
\(
\V\otimeshat\V\otimeshat \W \to
\W \otimeshat\V\otimeshat\V\).
\end{thm}
\begin{rem}
  In the same way as this section,
  we can construct factorized
  L-operators for $\Rcheck(\x)$ on $\Vm\otimeshat\Vm$
by using $V_n^-(\x)$ instead of $V_n(\x)$
  (cf.~Remark~\ref{rem:roperatorvm} and \ref{rem:incomingvm}).
In this case, outgoing intertwining vectors are characterized
by the following duality relation.
\begin{equation*}
\left\{
  \begin{aligned}
    &\sum_{i=k_1}^{k_2}\pk(\x)_{\l\,j}^{\l+\m\ve_i}
\pkbar(\x)_\l^{\l+\m\ve_i\,\ell}
=\d_{j\ell},\\
&\sum_{i=k_1}^{k_2}\pkbar(\x)_{\l}^{\l+\m\ve_j\,i}
\pk(\x)_{\l\,i}^{\l+\m\ve_\ell}
=\d_{j\ell}.
  \end{aligned}
\right.
\end{equation*}
Here,
for $\l,\k\in\L,\; k_1\leq j\leq k_2$,
we define $\pbar_{\bk}(\x)_\l^{\k\, j}\in\Bbb C$ as
follows (cf.~Remark~\ref{rem:roperatorvm}).
\[
  \pbar_{\bk}(\x)_\l^{\k\,j}:=
    \pbar_\l^\k(\vartheta {\frac12-\frac{j-k_1}k \brack{\frac{k}2}}
(\x+|\l|-kz,k\t)\exp (\p\rt-1 (k+1)z)).
\]
\end{rem}
\section{Vertex--irf correspondence and  factorized
  L-operators for Belavin's $R$-matrix}
In this section, we apply Theorem \ref{thm:vertex-irf} to
the $R$-matrix obtained through restricting the domain of
the elliptic $R$-operator to some finite-dimensional
subspace.
Then we will show that the vertex--IRF correspondence for
Belavin's $R$-matrix proved by Baxter \cite{baxter}, Jimbo,
Miwa and Okado
\cite{jimbo} is obtained from Theorem \ref{thm:vertex-irf}.
Moreover we will construct the factorized L-operators for
Belavin's $R$-matrix obtained by Hasegawa \cite{hasegawa1},
Quano and Fujii \cite{quano}.
First let us state the results proved by Felder and
Pasquier \cite{felder} more precisely.

For $k=1,2,\ldots$, let $\tVk (\x)$ be a space of entire
functions $f$ of one variable such that
\begin{equation*}
  \begin{array}{l}
f(z+1)=(-1)^kf(z),\\
f(z+\t)=(-1)^k\exp(-2\p\rt-1 (kz-\x+\frac{k\t}{2}))f(z).
  \end{array}
\end{equation*}
We note that $\tVk (\x)\subset \V$ if $k$ is even and that
$\tVk (\x)\subset \Vm$ if $k$ is odd.
In the same fashion as Theorem~\ref{thm:felder} and
Remark~\ref{rem:roperatorvm}, we obtain
\[
\Rcheck(\x_{12})(\tVk (\x_1)\otimes \tVk (\x_2+\m))
\subset \tVk (\x_2)\otimes \tVk (\x_1+\m).
\]
The space $\tVk (\x)$ is of $k$ dimensions and a basis is
given by
\begin{equation*}
\{ e_j(\x)(z)
:=\vartheta {\frac12-\frac{j}{k} \brack{\frac{k}{2}}}(\x-kz,k\t)
\}_{j\in \Bbb Z/k\Bbb Z}.
\end{equation*}
For $k=1,2,\ldots$,
define a translation operator $T_k(\x)$ on the space of
all holomorphic functions on $\Bbb C$
\cite{felder} by
\[
(T_k(\x)f)(z):=f(z-\frac{\x}{k}).
\]
$T_k(\x)$ maps isomorphically $\tVk:=\tVk (0)$
onto $\tVk (\x)$.
We modify the elliptic $R$-operator as
\[
\Rcheck_k(\x_{12}):=T_k(\x_2)^{-1}\otimes T_k(\x_1+\m)^{-1}
\Rcheck(\x_{12})T_k(\x_1)\otimes
T_k(\x_2+\m){\Bigg\vert}_{\tVk \otimes \tVk}.
\]
We note that $\Rcheck_k(\x_{12})$ is determined by the
difference $\x_{12}$.
In fact,
\begin{equation*}
\begin{split}
(\Rcheck_k(\x)f)(z_1,z_2)&=\frac{\vt1(\x)\vt1(z_{21}+\frac{\x
    +\m}{k}-\m)\vt1'(0)}{\vt1(-\m)\vt1(z_{21}+\frac{\x+\m}{k})}
f(z_2+\frac{\m}{k},z_1-\frac{\m}{k})\\
&\phantom{=}
+
\frac{\vt1(z_{21}+\frac{\x+\m}{k}-\x)\vt1'(0)}{\vt1(z_{21}
+\frac{\x+\m}{k})}f(z_1-\frac{\x}{k},z_2+\frac{\x}{k}).
\end{split}
\end{equation*}
Felder and Pasquier prove
\begin{thm}[\cite{felder}]\label{thm:belavin'sR}
  $\Rcheck_k(\x)$ preserves $\tVk \otimes \tVk$ and obeys the
  Yang--Baxter equation \eqref{intertwine:ybe}.
\end{thm}

Let
$\{ e^j\}_{j\in\Bbb Z /k\Bbb Z}\subset \tilde{V}_k^*$
be the dual basis of
$\{e_j:=e_j(0)\}\subset \tilde{V}_k$;
\[
e^i(e_j)=\d_{ij}.
\]
Now we define an operator $\Rcheck_k(\x)^*$ on $\tilde{V}_k^*\otimes
\tilde{V}_k^*$, the transpose of $\Rcheck_k(\x)$ on $\tVk
\otimes \tVk$.
\[
(\Rcheck_k(\x)^*e^\g\otimes e^\d)(e_\a\otimes e_\b)
:= (e^\d\otimes e^\g)(\Rcheck_k(\x)e_\b \otimes e_\a).
\]
\begin{prop}[cf.~\cite{felder}]
  The $R$-matrix $\Rcheck_k(\x)^*$ is Belavin's $R$-matrix
  up to constant.
\end{prop}
\begin{pf}
Let $A$ and $B$ be operators on the space of all holomorphic
functions on $\Bbb C$ as
\begin{equation*}
\begin{array}{l}
  (Af)(z)=-f(z+\frac1k),\\
  (Bf)(z)=-\exp (2\p\rt-1 (z+\frac{\t}{2k}))
    f(z+\frac{\t}{k}).
\end{array}
\end{equation*}
The space $\tVk$ is invariant under the actions of $A$ and
$B$.
In fact,
$A$ and $B$ are expressed
on $\tVk$ as
\begin{equation*}
  \begin{array}{l}
  Ae_j=e_j \exp \frac{2\p\rt-1 j}{k},\\
  Be_j=e_{j+1}.
  \end{array}
\end{equation*}
We define operators $A^*$ and $B^*$ on $\tilde{V}_k^*$
to be the transposes of $A$ and $B$ on $\tVk$, respectively;
\begin{equation*}
\begin{array}{l}
  A^*e^j=e^j \exp \frac{2\p\rt-1 j}{k},\\
  B^*e^j=e^{j-1}.
\end{array}
\end{equation*}
To prove this proposition, it is enough to show the following
\cite{belavin,hasegawa1,hasegawa2}.
\begin{enumerate}
\item $\Rcheck_k(\x)^*$ is an entire
$\mathrm{End}(\tilde{V}_k^*\otimes
\tilde{V}_k^*)$-valued function in $\x$.
\item \( \Rcheck_k(\x)^*x\otimes x=x\otimes x \Rcheck_k(\x)^*
  \quad x=A^*,B^*.\)
\item \( \Rcheck_k(\x+1)^*=(1\otimes A^*)^{-1}
\Rcheck_k(\x)^*(A^*\otimes
  1)\times (-1).\)
\item \( \Rcheck_k(\x+\t)^*=(1\otimes B^*)^{-1}
\Rcheck_k(\x)^*(B^*\otimes
  1)\times (-\exp 2\p\rt-1
  (\x+\frac{\t}{2}-\frac{\m}{k}))^{-1}.
\)
\item \( \Rcheck_k(0)^*=\vt1'(0)\,\mathrm{id}.\)
\end{enumerate}
The operator $\Rcheck_k(\x)$ on $\tVk \otimes \tVk$ has the
properties below, which imply the properties (2), (3), (4),
and (5) above, respectively.
\begin{enumerate}
\addtocounter{enumi}{1}
\item \( \Rcheck_k(\x)x\otimes x=x\otimes x \Rcheck_k(\x)
  \quad x=A,B.\)
\item \( \Rcheck_k(\x+1)=(1\otimes A)\Rcheck_k(\x)(A\otimes
  1)^{-1}\times (-1).\)
\item \( \Rcheck_k(\x+\t)=(1\otimes B)\Rcheck_k(\x)(B\otimes
  1)^{-1}\times (-\exp 2\p\rt-1
  (\x+\frac{\t}{2}-\frac{\m}{k}))^{-1}.
\)
\item \( \Rcheck_k(0)=\vt1'(0)\,\mathrm{id}.\)
\end{enumerate}
The proof is quite straightforward, so we omit it.

To prove (1), it suffices to show that
$\Rcheck_k(\x)$ is an entire $\mathrm{End}(\tVk \otimes
  \tVk)$-valued function in $\x$.
Let us introduce another basis of $\tVk$
(cf.~\cite{felder});
\[
\{ \tilde{e}_j(z):=
(-1)^j
\vartheta {\frac{k}{2} \brack{\frac12-\frac{j}{k}}}(z,\frac{\t}{k})
\}_{j \in \Bbb Z / k \Bbb Z}.
\]
In the same way as \cite{felder}, we can calculate the
matrix coefficients of
$\Rcheck_k(\x)$ on $\tVk \otimes \tVk$ with respect to the
basis $\{ \tilde{e}_i\otimes\tilde{e}_j \}$
and can check that all matrix coefficients are entire in $\x$.
This completes the proof.
\end{pf}

For $\l,\k\in\L$, we put
\(
\f(\x)_\l^\k:=\pbar_\k^\l \circ T_k(\x+|\l|_{\bk}-k\m){\Bigg\vert}_{\tVk}.
\)
Since
\begin{align*}
&\pbar_\k^\l \circ T_k(\x+|\l|_{\bk}-k\m)(e_j)\\
=&
\left\{
\begin{array}{ll}
\vartheta {\frac12-\frac{j}{k}
  \brack{\frac{k}{2}}}(\x+|\l|_{\bk}-k\l_i,k\t)
&
\text{if $\k=\l-\m\ve_i \ \ (k_1 \leq \exists i \leq k_2)$,}
\\
0,
&
\text{otherwise,}
\end{array}
\right.
\end{align*}
we get
\begin{align*}
\f(\x)_\l^\k
&=\sum_{j=0}^{k-1}
\pbar_\k^\l \circ T_k(\x+|\l|_{\bk}-k\m)(e_j)e^j\\
&=
\left\{
\begin{array}{ll}
\sum_{j=0}^{k-1}
\vartheta {\frac12-\frac{j}{k}
  \brack{\frac{k}{2}}}(\x+|\l|_{\bk}-k\l_i,k\t)
e^j,
&
\text{if $\k=\l-\m\ve_i \ \ (k_1 \leq \exists i \leq k_2)$,}
\\
0,
&
\text{otherwise.}
\end{array}
\right.
\end{align*}
Hence the vector $\f(\x)_\l^\k$ is nothing but the outgoing
intertwining vector of Belavin's $R$-matrix
\cite{hasegawa1,hasegawa2}, which was first
discovered by Baxter \cite{baxter}, Jimbo, Miwa and Okado \cite{jimbo}.

On the other hand, we put
\[
\tilde{W}
\begin{bmatrix}
  &\k&\\
\l&\x&\n\\
&\k'&
\end{bmatrix}
:=
\Wcheck
\begin{bmatrix}
  &\k'&\\
\n&\x&\l\\
&\k&
\end{bmatrix},
\]
and then
Theorem~\ref{thm:vertex-irf} and Remark~\ref{rem:incomingvm}
lead us to
\begin{thm}[Vertex--IRF correspondence for Belavin's
  $R$-matrix \cite{baxter,jimbo}]
For $\l,\k,\n\in\L$,
\begin{equation*}
\Rcheck_k(\x_{12})^*\f(\x_1)_\l^\k \otimes \f(\x_2)_\k^\n
=\sum_{\k'\in\L}
\f(\x_2)_\l^{\k'} \otimes \f(\x_1)_{\k'}^\n
\tilde{W}
\begin{bmatrix}
  &\k&\\
\l&\x_{12}&\n\\
&\k'&
\end{bmatrix}.
\end{equation*}
\end{thm}

Next we construct the factorized L-operators for Belavin's
$R$-matrix proved by Hasegawa
\cite{hasegawa1}, Quano and Fujii \cite{quano}.
To begin with, we introduce outgoing intertwining vectors in
$\tVk(\x)$ in the same fashion as Definition \ref{defn:outgoing}.
In the sequel, we fix $k_1, k_2 \in \Bbb Z$ such that
$k=k_2-k_1+1$ and assume that $\l,\k,\n \in \L$ and that $\x,\x_1,\x_2
\not\in \Bbb Z + \Bbb Z \t$.

For $k_1\leq j\leq k_2$,
we define $\vpkbar (\x)_\l^{\k\, j}\in\Bbb C$ by
\[
  \vpkbar (\x)_\l^{\k\,j}:=
    \pbar_\l^\k(e_j(\x+|\l|_{\bk})),
\]
and also define $\vpk(\x)_{\l\, j}^\k \in\Bbb C$ by the
following duality relations (cf.~Proposition
\ref{prop:weyl-kac});
\begin{equation*}
\left\{
  \begin{aligned}
    &\sum_{i=k_1}^{k_2}\vpk(\x)_{\l\,j}^{\l+\m\ve_i}
\vpkbar(\x)_\l^{\l+\m\ve_i\,\ell}
=\d_{j\ell},\\
&\sum_{i=k_1}^{k_2}\vpkbar(\x)_{\l}^{\l+\m\ve_j\,i}
\vpk(\x)_{\l\,i}^{\l+\m\ve_\ell}
=\d_{j\ell}.
  \end{aligned}
\right.
\end{equation*}
For $\k\ne\l+\m\ve_i\;(k_1\leq \forall i\leq k_2)$ we set
\[
  \vpk(\x)_{\l\ j}^\k:=0.
\]

Outgoing intertwining vectors
\(
\vpk(\x)_\l^\k(z)\in \tVk(\x+|\l|_{\bk})\)
of the elliptic $R$-operator are defined as
\[
  \vpk(\x)_\l^\k(z):=
    \sum_{j=k_1}^{k_2}\vpk(\x)_{\l\, j}^\k
e_j(\x+|\l|_{\bk})(z).
\]
Then we define the operators
\(
\Lcheckk(\x)_\l^\k
\)
as follows.
\[ (\Lcheckk(\x)_\l^\k f)(z):= \vpk(\x)_\l^\k(z)\pbar_\l^\k f,
\]
where $f \in \V$ if $k$ is even and $f \in \Vm$ if $k$ is odd.
In the same way as Section 3, these operators satisfy (cf.~Lemma
\ref{lem:loperator})
  \begin{equation*}
    \sum_{\k\in\L}\Rcheck(\x_{12})\Lcheckk(\x_1)_\l^\k
\otimes \Lcheckk(\x_2)_\k^\n
=\sum_{\k\in\L}
\Lcheckk(\x_2)_\l^\k \otimes \Lcheckk(\x_1)_\k^\n
\Rcheck(\x_{12}).
  \end{equation*}

We put
\[
\tilde{L}_{\bk}(\x)_\l^\k
:=
T_k(\x+|\k|_{\bk}-k\m)^{-1}
\Lcheckk(\x-k\m)_\k^\l
T_k(\x+|\k|_{\bk}-k\m)
{\Bigg\vert}_{\tVk}
\]
and denote its transpose as $\tilde{L}_{\bk}^*(\x)_\l^\k:
\tilde{V}_k^* \to \tilde{V}_k^*$.
Thus, for Belavin's $R$-matrix
$\Rcheck_k(\x)^*$,
  \begin{equation*}
    \sum_{\k\in\L}\Rcheck_k(\x_{12})^* \tilde{L}_{\bk}^*(\x_1)_\l^\k
\otimes \tilde{L}_{\bk}^*(\x_2)_\k^\n
=\sum_{\k\in\L}
\tilde{L}_{\bk}^*(\x_2)_\l^\k \otimes \tilde{L}_{\bk}^*(\x_1)_\k^\n
\Rcheck_k(\x_{12})^*.
  \end{equation*}
We define an operator $\tilde{L}_{\bk}^*(\x):\tilde{V}_k^* \otimes
\W \to \W \otimes \tilde{V}_k^*$ by
\begin{equation*}
\tilde{L}_{\bk}^*(\x)(e^i\otimes \d^\k)
= \sum_{\l\in\L}\d^\l\otimes \tilde{L}_{\bk}^*(\x)_\l^\k e^i.
\end{equation*}
The theorem below tells us that the operator
$\tilde{L}_{\bk}^*(\x)$ is
the factorized L-operator for Belavin's $R$-matrix,
which were first constructed by
Hasegawa \cite{hasegawa1}, Quano and Fujii \cite{quano}.
\begin{thm}[Factorized L-operator for Belavin's $R$-matrix]
For $\x_1,\x_2\not\in\Bbb Z + \Bbb Z \t$
  \begin{equation*}
(1\otimes \Rcheck_k(\x_{12})^*)(\tilde{L}_{\bk}^*(\x_{1})\otimes
1)(1\otimes \tilde{L}_{\bk}^*(\x_{2}))
=(\tilde{L}_{\bk}^*(\x_{2})\otimes 1)(1\otimes
\tilde{L}_{\bk}^*(\x_{1}))(\Rcheck_k(\x_{12})^* \otimes 1).
  \end{equation*}
Here both sides are the operators \(
\tilde{V}_k^* \otimes \tilde{V}_k^*\otimes \W \to
\W \otimes \tilde{V}_k^* \otimes \tilde{V}_k^* \).
\end{thm}
\begin{ack}
  The author would like to express his deep gratitude
to Professor Kimio Ueno for useful advice,
to Professor Takayuki Hibi for
constant encouragement, and
to Professor Michio Jimbo for pointing him out another proof
of Corollary \ref{cor:star-triangle} which is used in this
paper.
\end{ack}
\end{document}